\documentclass[12pt]{amsart}

\usepackage{amssymb,latexsym,amscd}
\usepackage{amsthm,amsmath}

\newcommand{\gfr}{\mathfrak{g}}
\newcommand{\hfr}{\mathfrak{h}}
\newcommand{\pfr}{\mathfrak{p}}

\newcommand{\cV}{\mathcal{V}}
\newcommand{\cH}{\mathcal{H}}

\newcommand{\cO}{\mathcal{O}}

\newcommand{\cG}{\mathcal{G}}

\newcommand{\cM}{\mathcal{M}}

\newcommand{\cX}{\mathcal{X}}
\newcommand{\cL}{\mathcal{L}}

\newcommand{\Vir}{\mathrm{Vir}}
\newcommand{\mult}{\mathrm{mult}}
\newcommand{\Spec}{\mathrm{Spec}}

\newcommand{\lra}{\longrightarrow}

\newcommand{\ra}{\rightarrow}

\newcommand{\PP}{\mathbb{P}}

\newcommand{\ZZ}{\mathbb{Z}}

\newcommand{\CC}{\mathbb{C}}

\newcommand{\End}{\mathrm{End}}

\newcommand{\SL}{\mathrm{SL}}

\newcommand{\Spin}{\mathrm{Spin}}

\newcommand{\Lie}{\mathrm{Lie}}

\theoremstyle{plain}

\newtheorem{thm}{Theorem}[section]

\newtheorem{lem}[thm]{Lemma}

\newtheorem{prop}[thm]{Proposition}

\newtheorem{rem}[thm]{Remark}
\newtheorem{defi}[thm]{Definition}

\begin{document}

\title[]{Strange duality for Verlinde spaces of exceptional groups at level one}

\begin{abstract}
The moduli stack $\cM_X(E_8)$ of principal $E_8$-bundles over a smooth projective curve $X$ carries a
natural divisor $\Delta$. We study the pull-back of the divisor $\Delta$ to the moduli stack
$\cM_X(P)$, where $P$ is a semi-simple and simply connected group such that its Lie algebra $\mathrm{Lie}(P)$ is a
maximal conformal subalgebra of $\mathrm{Lie}(E_8)$. We show that the divisor $\Delta$ induces
``Strange Duality"-type isomorphisms between the Verlinde spaces at level one of the following pairs of groups
$(\SL(5), \SL(5))$, $(\Spin(8), \Spin(8))$, $(\SL(3), E_6)$ and $(\SL(2), E_7)$.

\end{abstract}

\author{Arzu Boysal}

\author{Christian Pauly}

\address{Matematik B\"{o}l\"{u}m\"{u}\\ Bo\={g}azi\c{c}i
\"{U}niversitesi\\
TR-34342 \\ Bebek \\ Istanbul \\
T\"{u}rkiye (Turkey).}

\email{arzu.boysal@boun.edu.tr}

\address{D\'epartement de Math\'ematiques \\ Universit\'e de Montpellier II - Case Courrier 051 \\ Place Eug\`ene Bataillon \\ 34095 Montpellier Cedex 5 \\ France.}

\email{pauly@math.univ-montp2.fr}




\maketitle

\section{Introduction}
Let $X$ be a smooth complex projective curve of genus $g$ and let $G$ be a simple and simply connected
complex Lie group. We denote by $\cM_X(G)$ the moduli stack parametrizing principal $G$-bundles over the
curve $X$ and by $\cL_G$ the ample line bundle over $\cM_X(G)$ generating its Picard group. The starting point of our
investigation is the observation (see e.g. \cite{So}, \cite{F1}, \cite{F2}) that
$$ \dim H^0(\cM_X(E_8), \cL_{E_8}) = 1. $$
for any genus $g$. In other words, the moduli stack $\cM_X(E_8)$ carries a natural
divisor $\Delta$. Unfortunately a
geometric interpretation of this divisor is not known.

\bigskip

In this paper we study the pull-back of this mysterious divisor $\Delta$ under the morphisms
$\cM_X(P) \rightarrow \cM_X(E_8)$ induced by the group homomorphisms $\phi: P \rightarrow E_8$, where we assume that
$P$ is connected, simply connected and semi-simple, and that the differential $d\phi: \pfr = \Lie(P)  \rightarrow \mathfrak{e_8} = \Lie(E_8)$
is a {\em conformal} embedding of Lie algebras (see Definition~\ref{cb}). 
We recall (\cite{BB} p.~566) that any subalgebra of maximal rank $8$ of $\mathfrak{e_8}$ (see \cite{BD} Chapter $7$
for a list) is actually a conformal subalgebra of $\mathfrak{e_8}$ with Dynkin (multi-)index one.  Maximal conformal
subalgebras of $\mathfrak{e_8}$ with Dynkin (multi-)indicies one have been classified by \cite{BB} and \cite{SW},
and the full list is as follows:
\begin{equation} \label{confE8}
\begin{array}{ll}
\text{maximal rank} &: \mathfrak{so}(16), \ \mathfrak{sl}(9), \ \mathfrak{sl}(5) \oplus \mathfrak{sl}(5), \ \mathfrak{sl}(3) \oplus
\mathfrak{e_6}, \  \mathfrak{sl}(2) \oplus \mathfrak{e_7}\\
\text{non-maximal rank}&: \ \mathfrak{g}_2 \oplus \mathfrak{f}_4.
\end{array}
\end{equation}
In Table (\ref{listP}) we list the corresponding simply connected
Lie groups $P$ and the finite kernel $N$ of their natural maps to
$E_8$ (see e.g. \cite{CG} Lemma 3.3).

\begin{equation} \label{listP}
\begin{tabular}{|c||c|c|c|c|c|c|}
        \hline
            P & $\Spin(16)$ &  $\SL(9)$ &  $\SL(5) \times \SL(5)$ & $\SL(3) \times E_6$ &  $\SL(2) \times E_7$ & $G_2 \times F_4$ \\
            \hline
            N & $\ZZ/2 \ZZ$ &  $\ZZ/ 3\ZZ$ & $\ZZ/5\ZZ$ & $\ZZ/3 \ZZ$ & $\ZZ/ 2 \ZZ$ & 1 \\
     \hline

\end{tabular}
\end{equation}

\bigskip

Note that $N$ is a subgroup of the center of $P$.
We introduce the finite abelian group $\cM_X(N)$ of principal $N$-bundles over $X$, which acts
on $\cM_X(P)$ by twisting $P$-bundles with $N$-bundles. Since $N$ is the kernel of $\phi$, the
group $\cM_X(N)$ acts on the fibers of the induced stack morphism $\tilde{\phi} : \cM_X(P) \rightarrow \cM_X(E_8)$.
We emphasize that the $\cM_X(N)$-linearization of the line bundle $\cL_P$ over $\cM_X(P)$ is not unique, but  we obtain a 
{\em canonical} $\cM_X(N)$-linearization (hence a {\em canonical} linear action
of $\cM_X(N)$ on $H^0(\cM_X(P), \cL_P)$), namely, the one induced from the isomorphism $\tilde{\phi}^* \cL_{E_8} = \cL_P$.

If $P$ has two simple factors, the line bundle
$\cL_P$ denotes the tensor product of the ample generators on each factor (see section 2.1).
With this notation our main result can be stated as follows.

\bigskip

\noindent
{\bf Theorem 1.} Let $P$ be a Lie group of  Table \eqref{listP} and let
$$\phi_P : H^0(\cM_X(E_8), \cL_{E_8}) \longrightarrow H^0(\cM_X(P), \cL_P)$$
be the linear map on global sections induced by the homomorphism $\phi: P \rightarrow E_8$. Then
\begin{itemize}
\item[(i)] The linear map $\phi_P$ is nonzero.
\item[(ii)] In all cases except $P = G_2 \times F_4$, its one-dimensional image coincides with the
$\cM_X(N)$-invariant subspace of $H^0(\cM_X(P), \cL_P)$, where $\cM_X(N)$ acts on $H^0(\cM_X(P), \cL_P)$
with the canonical action.
\end{itemize}

\bigskip

Unfortunately, even in the cases $\SL(9)$ and $\Spin(16)$, we are not able to give a geometric description of the
zero-divisor of the $\cM_X(N)$-invariant section in the moduli stacks $\cM_X(\SL(9))$ and $\cM_X(\Spin(16))$ ---
see section 7.1 for further discussion.

\bigskip

A word about the proof of Theorem 1. We make use of the
identification of the space of generalized $G$-theta functions
$H^0(\cM_X(G), \cL_G)$ with the space of conformal blocks
$\cV_0^\dagger(X,\gfr)$ associated to the curve $X$ with one marked
point labelled with the zero weight. Here $\gfr = \mathrm{Lie}(G)$.
Under this identification the linear map $\phi_P$ becomes the map
induced by the natural inclusion of the basic highest weight modules
$\cH_0(\pfr) \hookrightarrow \cH_0(\mathfrak{e}_8)$ of the affine
Lie algebras $\hat{\pfr}$ and $\hat{\mathfrak{e}_8}$. The proof has
essentially two steps. First, we use a result by P. Belkale
(\cite{B} Proposition 5.8) saying that the linear map $\phi_P$ has
constant rank when the curve $X$ varies in a family of smooth
curves. Here, the fact that the embedding $\mathfrak{p} \subset
\mathfrak{e}_8$ is conformal, is crucial, since it ensures that
$\phi_P$ is projectively flat with respect to the WZW connections on
both sheaves of vacua over any family of smooth curves. Secondly, we
study in section 4 the behaviour of the factorization rules of the
spaces of conformal blocks associated to $\gfr$ under any conformal
embedding $\mathfrak{p} \subset \gfr$. This will follow once one has
decomposed the ``sewing-procedure" tensor $\tilde{\gamma}_\lambda
\in \cH_\lambda(\gfr) \otimes \cH_{\lambda^\dagger}(\gfr)[[q]]$
under the decomposition of the $\hat{\gfr}$-modules
$\cH_\lambda(\gfr)$ and $\cH_{\lambda^\dagger}(\gfr)$ into
irreducible $\hat{\pfr}$-modules. Finally, combining these two steps
allows us to show by induction on the genus of $X$ that $\phi_P$ is
non-zero.

\bigskip

In the cases when $P$ is not simple and $N$ not trivial, an argument using the representation theory of
Heisenberg groups allows us to show the following result, which can be seen as an instance of Strange Duality
for exceptional groups at level one.

\bigskip
\noindent
{\bf Theorem 2.} Let $(A,B)$ be one of the three pairs $(\SL(5), \SL(5)),(\SL(3),E_6), (\SL(2),E_7)$.  Consider {\em any}
$\cM_X(N)$-linearization of the line bundle $\cL_{A \times B}$ over $\cM_X(A \times B)$.  Let
$$\sigma \in H^0(\cM_X(A \times B), \cL_{A \times B})= H^0(\cM_X(A), \cL_A) \otimes H^0(\cM_X(B), \cL_B)$$ be a non-zero
element of the one-dimensional $\cM_X(N)$-invariant subspace. Then  $\sigma$ induces an isomorphism
$$ \sigma : H^0(\cM_X(A), \cL_A)^* \longrightarrow H^0(\cM_X(B), \cL_B).$$

\bigskip

A similar isomorphism is obtained for the pair $(\Spin(8),
\Spin(8))$ --- see section 7.2.1. 

We would like to mention that
the proofs of Theorem 1 and Theorem 2 are independent, and that both results are
related by the fact that the ``Strange Duality" isomorphism of Theorem 2 corresponding 
to the {\em canonical} $\cM_X(N)$-linearization is obtained precisely by pulling-back
the $E_8$-theta divisor $\Delta$ to the moduli stack $\cM_X(A \times B)$.

\bigskip

The last few years have seen important progress on ``Strange Duality" or ``rank-level" duality
for Verlinde spaces. For a survey we refer, for example, to the papers \cite{MO}, \cite{Po} or \cite{Pa}.

\bigskip

We would like to thank Laurent Manivel and Nicolas Ressayre for helpful comments, as well as the referees for helping us 
improve readability of the paper.
\bigskip

\noindent
{\bf Acknowledgements.} The first author wishes to thank the Consejo Superior de Investigaciones Cient\'ificas (CSIC, Spain) and the Universit\'e de Montpellier II (France) for financial support of research visits.
The second author was partially supported by the Ministerio de Educaci\'on y Ciencia (Spain) through the grant
SAB2006-0022.

\section{Notation and preliminaries}

\subsection{Moduli stacks and line bundles}

\subsubsection{Dynkin index}

Let $\pfr$ and $\gfr$ be two simple Lie algebras and let $\varphi: \pfr \ra \gfr$ be a Lie algebra
homomorphism. There exists \cite{D} a unique integer $d_\varphi$, called the Dynkin index of the
homomorphism $\varphi$, satisfying
$$ (\varphi(x),\varphi(y))_\gfr = d_\varphi (x,y)_\pfr, \qquad \text{for all} \ x,y \in \pfr,$$
where $( \ , \ )_*$ denotes the Cartan-Killing form on $\pfr$ and $\gfr$, normalized such that
$( \theta, \theta ) = 2$ for their respective highest roots $\theta$. If $\pfr$ is
semi-simple with two components $\pfr_1 \oplus \pfr_2$, then the Dynkin multi-index of $\varphi =
\varphi_1 \oplus \varphi_2 : \pfr_1 \oplus \pfr_2 \ra \gfr$ is given by $d_\varphi = (d_{\varphi_1},
d_{\varphi_2})$, where each $d_{\varphi_i}$ is defined using $\varphi_i : \pfr_i \ra \gfr$.

\subsubsection{Line bundles over the moduli stack $\cM_X(P)$}

If $P$ is a simple and simply connected complex Lie group we refer to \cite{LS} and \cite{So} for the description
of the ample generator $\cL_P$ of the Picard group of the moduli stack $\cM_X(P)$. If $P = P_1 \times
P_2$ with $P_i$ simple and simply connected, we put $\cL_P = \cL_{P_1} \boxtimes \cL_{P_2}$ and we note
$\cL_P^{d} = \cL_{P_1}^{d_1} \boxtimes \cL_{P_2}^{d_2}$ for a multi-index $d=(d_1,d_2)$. The following
lemma follows easily from \cite{LS} and \cite{KNR}.

\begin{lem}
Let $\phi: P \ra G$ be a homomorphism between simply-connected complex Lie groups with $G$ simple and $P$ semi-simple. Let $\tilde{\phi} : \cM_X(P) \ra \cM_X(G)$ be the induced stack morphism. Then we
have the equality
$$ \tilde{\phi}^* \cL_G = \cL_P^{d_\varphi},$$
where $d_\varphi$ is the Dynkin (multi-) index of the differential $\varphi = d\phi: \pfr \ra \gfr$.
\end{lem}

\subsection{Spaces of conformal blocks}

\subsubsection{The case $\gfr$ simple}

Let $\gfr$ be a simple Lie algebra and $\hfr \subset \gfr$ a Cartan subalgebra. We denote 
as before by $( \ , \ )$
the normalized Cartan-Killing form and we will use the same notation for the restricted form on
$\hfr$ and for the induced form on $\hfr^*$. We consider (see \cite{K} Chapter 7) the non-twisted affine Lie algebra associated to $\gfr$
over $\CC ((z))$
$$\hat{\gfr} = \gfr \otimes \CC ((z)) \oplus \CC c \oplus \CC d$$
with Lie bracket
$$[x\otimes f,y\otimes g] = [x,y]\otimes fg + (x,y) \mathrm{Res}_{z=0} (gdf) \cdot c, \ [\gfr,c] = 0, \
[d,c] = 0, \ [d,x(n)] = nx(n),$$
for $x,y \in \gfr, \ f,g \in \CC((z))$ and $ n \in \ZZ$. Here we put $x(n) = x \otimes z^n$. We identify
$\gfr$ with the subalgebra $\gfr \otimes 1$ of $\hat{\gfr}$. The subalgebra
$$\hat{\hfr} = \hfr \oplus \CC c \oplus \CC d$$
is the Cartan subalgebra of the affine Lie algebra $\hat{\gfr}$. We extend $\lambda \in \hfr^*$ to
a linear form on $\hat{\hfr}$ by putting $\langle \lambda, \CC c \oplus \CC d \rangle = 0$, where
$\langle \ , \ \rangle$ is the standard pairing. We
define the elements $\Lambda_0$ and $\delta$ in the dual $\hat{\hfr}^* = \hfr^* \oplus \CC \Lambda_0
\oplus \CC \delta$ by $\langle \delta , d\rangle = \langle \Lambda_0, c \rangle = 1$ and
$\langle \delta, \hfr \oplus \CC c \rangle = \langle \Lambda_0, \hfr \oplus \CC d \rangle = 0$.
We extend the form $( \ ,
\ )$ to $\hat{\hfr}^*$ by putting
$$ (\hfr^* , \CC \Lambda_0 \oplus \CC \delta) = 0, \qquad (\delta,\delta) = (\Lambda_0, \Lambda_0) = 0,
\qquad (\delta,\Lambda_0) = 1.$$

\bigskip

The Weyl group of $\gfr$ is denoted by $W(\gfr)$. Call $w_0^\gfr \in W(\gfr)$ its longest element. Later we will need the following fact

\begin{prop} \label{Weylgroup}
Let $\pfr \subset \gfr$ be an embedding of semi-simple Lie algebras and choose Cartan subalgebras
such that $\hfr_\pfr \subset \hfr_\gfr$. Then there exists an element $\tilde{w} \in W(\gfr)$ which preserves the
subspace $\hfr_\pfr \subset \hfr_\gfr$ and such that the restriction $\tilde{w}_{|{\hfr_\pfr}}$ coincides with the
longest element $w_0^\pfr$.
\end{prop}

\begin{proof}
This can be deduced from a more general fact --- see e.g. Theorem 2.1.4 \cite{BS}. Moreover if $\pfr$ and
$\gfr$ have the same rank (all our cases except $\mathfrak{g}_2 \oplus \mathfrak{f}_4$), i.e.  $\hfr_\pfr = \hfr_\gfr$, we have a canonical inclusion $W(\pfr) \subset W(\gfr)$.
\end{proof}

\bigskip

Next, we recall some representation theory of the affine Lie
algebra $\hat{\gfr}$ from \cite{K} Chapter 12. 

\bigskip

We denote by $P(\gfr)$ the weight lattice of $\gfr$ and by $P_+(\gfr)$ the subset of $P(\gfr)$ consisting of dominant integral weights.
Given a positive integer $k$, called level,  we consider the finite set
$$P_k(\gfr) := \{ \lambda \in P_+(\gfr) \ | \ ( \lambda , \theta )  \leq k \} \subset \hfr^*. $$
Given $\lambda \in P_k(\gfr)$ we denote by 
$\cH_\lambda(\gfr)$, or simply $\cH_\lambda$ if no confusion arises, the integrable
$\hat{\gfr}$-module with highest weight $\lambda+k\Lambda_0$; in particular,
\medskip

(1) the center $c \in \hat{\hfr}$ acts on $\cH_\lambda$ as $k \cdot \mathrm{Id}$. 

(2) the derivation $d \in \hat{\hfr}$ acts trivially on the highest weight vector $v_{\lambda}$ of $\cH_\lambda$.

\bigskip

We introduce the set 
$$ \widehat{P_k(\gfr)} := \{ \hat{\lambda} = \lambda + k \Lambda_0 + \zeta \delta \ | \ \lambda \in P(\gfr) ,
\zeta \in \CC \} \subset \hat{\hfr}^*.$$
Given $\hat{\lambda}=\lambda+k\Lambda_0+\zeta \delta$ in $\widehat{P_k(\gfr)}$ we define 
$\hat{\lambda}^\dagger=-w_0^\gfr(\lambda)+k\Lambda_0+\zeta \delta$.  This gives an involution
$\hat{\lambda}\to \hat{\lambda}^\dagger$ on the set $\widehat{P_k(\gfr)}$.

Note that there is a projection map
$$\widehat{P_k(\gfr)} \rightarrow P(\gfr), \qquad  \hat{\lambda} = \lambda + k \Lambda_0 + \zeta \delta \mapsto \lambda.$$
We will view $P_k(\gfr)$ as a subset of $\widehat{P_k(\gfr)}$ under the mapping $\lambda \mapsto
\hat{\lambda} = \lambda + k \Lambda_0$ and we observe that the involution $\hat{\lambda}\to \hat{\lambda}^\dagger$ restricts
to an involution $\lambda \mapsto
\lambda^{\dagger}= -w_0^\gfr(\lambda)$ on the finite set $P_k(\gfr)$.  Note that $-\lambda^\dagger$ is the lowest weight of
the irreducible right $\gfr$-module $V_\lambda^*$, the dual of
$V_\lambda$. 

\bigskip

More generally, for any $\hat{\lambda} = \lambda + k \Lambda_0 + \zeta \delta \in \widehat{P_k(\gfr)}$ such that $\lambda
\in P_k(\gfr)$ we denote by $\cH_{\hat{\lambda}}(\gfr)$ the integrable  $\hat{\gfr}$-module with highest weight 
$\hat{\lambda} = \lambda + k \Lambda_0 + \zeta \delta$; in particular,

\medskip

(1) the center $c \in \hat{\hfr}$ acts on $\cH_\lambda$ as $k \cdot \mathrm{Id}$. 

(2) the derivation $d \in \hat{\hfr}$ acts on the highest weight vector $v_{\lambda}$ of $\cH_\lambda$ as 
$d\cdot v_{\lambda}=\zeta v_{\lambda}$.

\bigskip

For $\hat{\lambda} = \lambda + k \Lambda_0 + \zeta \delta \in \widehat{P_k(\gfr)}$ with $\lambda
\in P_k(\gfr)$, we note that the $\hat{\gfr}$-modules $\cH_\lambda$ and
$\cH_{\hat{\lambda}}$ become  isomorphic as modules over the subalgebra
$[\hat{\gfr}, \hat{\gfr}] = \gfr \otimes \CC((z)) \oplus \CC c$. The $\hat{\gfr}$-module
$\cH_0(\gfr)$ with zero weight and level $1$ (i.e., $\hat{\lambda} = \Lambda_0$) is called the {\em basic} $\hat{\gfr}$-module.

\bigskip

Given $s$ points $\vec{p} = (p_1, \ldots, p_s)$ on $X$ we consider the open subset 
$U=X\setminus \{p_1, \ldots, p_s \}$, and choose a local coordinate $\xi_i$ at each
point $p_i$. 
Following \cite{U} Definition 3.1.1 we introduce the Lie algebra 
$$\hat{\gfr}_s:=\bigoplus_{j=1}^s \gfr\otimes_{\CC}\CC((\xi_i))\oplus \CC c .$$

There is a natural embedding of the ring $H^0(U, \mathcal{O}_U)$ of regular functions on $U$
into $\oplus_{j=1}^s \CC((\xi_i))$  
via Laurent expansions at the points $p_i$.   By \cite{U} Lemma 3.1.2,
$\gfr(U):=\gfr\otimes H^0(U, \mathcal{O}_U)$ is a Lie subalgebra of $\hat{\gfr}_s$.

For $\vec{\lambda} = (\lambda_1, \lambda_2, \ldots , \lambda_s) \in P_k(\gfr)^s$ we introduce the left $\hat{\gfr}_s$-module
$$\cH_{\vec{\lambda}} := \cH_{\lambda_1} \otimes \cH_{\lambda_2} \otimes \ldots \otimes \cH_{\lambda_s}. $$
The space of covacua
associated to the data $(X, \vec{p} , \vec{\lambda}, \vec{\xi})$ is defined as
$$\cV_{\vec{\lambda}}(X, \gfr) := \cH_{\vec{\lambda}} / \gfr(U) \cdot \cH_{\vec{\lambda}},$$
where $\gfr(U)$ acts on $\cH_{\vec{\lambda}}$ via the inclusion in $\hat{\gfr}_s$.
We refer to \cite{U} for further details.

The space of vacua or the {\em space of conformal blocks} is defined as the dual of $\cV_{\vec{\lambda}}(X, \gfr)$ and is denoted by
$\cV^{\dagger}_{\vec{\lambda}}(X,\gfr)$. We note that there is an inclusion
$$ \cV^{\dagger}_{\vec{\lambda}}(X,\gfr) \hookrightarrow \cH^\dagger_{\vec{\lambda}}.$$

\medskip

We denote by $\langle \ | \ \rangle$ the natural pairing between $\cH_{\vec{\lambda}}$ and its dual
$\cH_{\vec{\lambda}}^\dagger$.

\bigskip

The construction of the space of conformal blocks can be carried out for a family $\cX \ra S$ of pointed
nodal curves and provides a sheaf
over the base scheme $S$, called the sheaf of vacua and denoted  by $\cV^{\dagger}_{\vec{\lambda}}(\cX,\gfr)$.
A fundamental property is that the sheaf of vacua is locally free and that it commutes with any base change in S 
(see e.g. \cite{U} Theorem 4.4.2).

\subsubsection{The case $\gfr$ semi-simple}

We now adapt the previous contructions to semi-simple Lie algebras $\gfr$. For our purposes it is enough to deal with the case
when $\gfr$ is the direct sum of two simple Lie algebras $\gfr = \gfr_1 \oplus \gfr_2$.
By \cite{K} section 12.9 we define the affine Lie algebra associated to $\gfr$ by
$$ \hat{\gfr} = \gfr \otimes \CC ((z)) \oplus \CC c_1 \oplus \CC c_2 \oplus \CC d,$$
with Cartan subalgebra $\hat{\hfr} = \hfr \oplus \CC c_1 \oplus \CC c_2 \oplus \CC d$. Similar to the
case of simple algebras, one defines a Lie bracket and a non-degenerate bilinear form on $\hat{\hfr}$ and
on its dual $\hat{\hfr}^* = \hfr^* \oplus \CC \Lambda_0^{(1)} \oplus \CC \Lambda_0^{(2)} \oplus \CC \delta$.

\bigskip

Given a multi-index $k = (k_1, k_2)$,  with $k_i$ positive integers,
we introduce the sets
\begin{eqnarray*}
 P_{k}(\gfr) & = & P_{k_1} (\gfr_1) \times  P_{k_2} (\gfr_2),  \\
\widehat{P_k(\gfr)} & = & \{ \hat{\lambda} = \lambda + k_1
\Lambda_0^{(1)} + k_2 \Lambda_0^{(2)} + \zeta \delta \ | \ \lambda
\in P(\gfr) , \zeta \in \CC \} \subset \hat{\hfr}^*,
\end{eqnarray*}
and we associate to a weight $\lambda = (\lambda^{(1)}, \lambda^{(2)}) \in P_{k}(\gfr)$ the integrable
$\hat{\gfr}$-module
$$\cH_\lambda(\gfr)  = \cH_{\lambda^{(1)}}(\gfr_1) \otimes
\cH_{\lambda^{(2)}}(\gfr_2).$$
Similarly, we introduce for any $\hat{\lambda} \in \widehat{P_k(\gfr)}$ the $\hat{\gfr}$-module $\cH_{\hat{\lambda}}(\gfr)$.

With these definitions it is easy to deduce the following decomposition
of the spaces of conformal blocks
$$  \cV^{\dagger}_{\vec{\lambda}}(\gfr) =  \cV^{\dagger}_{\vec{\lambda}^{(1)}}(\gfr_1) \otimes  \cV^{\dagger}_{\vec{\lambda}^{(2)}}(\gfr_2) \qquad \text{with} \ \vec{\lambda} = (\vec{\lambda}^{(1)},\vec{\lambda}^{(2)}) \in P_{k}(\gfr)^s.$$

\bigskip

\subsection{Generalized theta functions and the Verlinde formula}

For the convenience of the reader we recall (see \cite{LS}, \cite{F3}, \cite{KNR})
that there is an isomorphism between the space
$H^0(\cM_X(G), \cL_G^k)$ of generalized $G$-theta functions of level $k$ and
the space of conformal blocks $\cV^{\dagger}_0(X, \gfr)$ associated to
the curve $X$ with one marked point labelled with the zero weight at level $k$, i.e,
\begin{equation} \label{verlindeiso}
H^0(\cM_X(G), \cL_G^k) \stackrel{\sim}{\longrightarrow} \cV^{\dagger}_0(X, \gfr).
\end{equation}
The dimension of this space is given by the Verlinde formula (see \cite{F3}, \cite{T}). Its value for the
groups $G$ of type $ADE$ and at level one is $|Z|^g$, where $|Z|$ denotes the order of the center $Z$ of $G$
(see \cite{F1}).  The list is as follows:

\begin{equation} \label{tableVerlinde}
\begin{tabular}{|c||c|c|c|c|c|}
        \hline
            $G $&  $\SL(n)$ & $\Spin(2n)$ &  $E_6$ &  $E_7$  & $E_8$  \\
            \hline
            $\dim H^0(\cM_X(G),\cL_G)$ & $n^g$ &  $4^g$ & $3^g$ & $2^g$ & $1$  \\
     \hline
     \end{tabular}
\end{equation}

\bigskip

\section{Conformal pairs: properties of their representations}

\subsection{The Virasoro algebra and its representation on $\cH_{\hat{\lambda}}(\gfr)$}

The Virasoro algebra
$$\Vir := \bigoplus_{n \in \ZZ} \CC d_n \oplus \CC \tilde{c},$$
is defined by the relations $[d_i,d_j] = (i-j)d_{i+j} + \frac{1}{12}(i^3-i)\delta_{i,-j} \tilde{c}$
and $[d_i,\tilde{c}]=0$.

\bigskip

If $\gfr$ is simple, we define for any level $k$ and any
$\hat{\lambda} \in \widehat{P_k(\gfr)}$, the conformal anomaly
$c(\gfr,k)$ and  the trace anomaly $ \Delta_{\hat{\lambda}}(\gfr)$
as
$$c(\gfr,k) = \frac{k \dim \gfr}{\check{h}(\gfr) + k}, \qquad \text{and} \qquad \Delta_{\hat{\lambda}}(\gfr) = \frac{(\hat{\lambda}, \hat{\lambda} + 2 \hat{\rho})}{2(\check{h}(\gfr) + k)}.$$
Here $\check{h}(\gfr)$ is the dual Coxeter number of $\gfr$ and $\hat{\rho} = \rho + \check{h}(\gfr) \Lambda_0$, where $\rho$ denotes 
the half-sum of the positive roots of $\gfr$. We choose dual bases $\{ u_i \}$ and $\{ u^i \}$ of the simple algebra $\gfr$ and introduce for any $n \in \ZZ$ the
Sugawara operator (see \cite{KW} section 3.2)
\begin{equation} \label{defL0}
L_{n}^\gfr = \frac{1}{2(k+ \check{h}(\gfr))} \sum_{j \in \ZZ}
\sum_i : u_i(-j) u^i(j+n) : ,
\end{equation}
where the notation $ : \ \ :$ stands for the normal ordering. Then $L_{n}^\gfr$ acts linearly on
$\cH_{\hat{\lambda}}(\gfr)$ and, by putting $d_n = L_{n}^\gfr$ and
$\tilde{c} = c(\gfr,k) \mathrm{Id}$, we obtain a representation of
$\Vir$ on $\cH_{\hat{\lambda}}(\gfr)$.
\bigskip

If $\gfr$ is semi-simple with $\gfr = \gfr_1 \oplus \gfr_2$, we define for $k = (k_1,k_2)$, see \cite{KW} formula (1.4.7),
$$c(\gfr,k) = c(\gfr_1,k_1) + c(\gfr_2,k_2), \qquad \text{and} \qquad \Delta_{\hat{\lambda}}(\gfr) =
\Delta_{\hat{\lambda}^{(1)}}(\gfr_1) +
\Delta_{\hat{\lambda}^{(2)}}(\gfr_2),$$ and we put $L_{n}^\gfr =
L_{n}^{\gfr_1} + L_{n}^{\gfr_2}$. As in the simple case, we obtain a
representation of $\Vir$ on $\cH_{\hat{\lambda}}(\gfr)$.

\bigskip

For later use we recall  the following relation (\cite{KW} formula (1.4.5))
\begin{equation} \label{traceanomaly}
 \Delta_{\hat{\lambda} + n \delta}(\gfr) = \Delta_{\hat{\lambda}}(\gfr) + n .
\end{equation}

\bigskip

The endomorphism $L_0^\gfr$ of the $\hat{\gfr}$-module $\cH_{\hat{\lambda}}(\gfr)$ can be diagonalized
$$  \cH_{\hat{\lambda}} = \bigoplus_{m=0}^{\infty} \cH_{\hat{\lambda}}(m) \qquad \text{with} \qquad
\cH_{\hat{\lambda}}(m) := \{
u \in \cH_{\hat{\lambda}} \ | \ L_0^\gfr(u) = (\Delta_{\hat{\lambda}}(\gfr) + m) u \}.$$

We recall that the endomorphism $L_0^\gfr$ of $\cH_{\hat{\lambda}}(\gfr)$ defined by \eqref{defL0} can also be written as   (\cite{KW} (3.2.6) or \cite{K} Corollary 12.8 (b))
\begin{equation} \label{L0}
L_0^\gfr = \Delta_{\hat{\lambda}}(\gfr) \mathrm{Id} - d.
\end{equation}
Note that $\cH_{\hat{\lambda}}(0)$ equals the irreducible $\gfr$-module $V_\lambda$ with highest
weight $\lambda$.

\subsection{Definition of conformal pair}

\begin{defi}[e.~g.~\cite{K} Chapter 13] \label{cb}
Let $\pfr$ be a semisimple subalgebra of a simple Lie algebra $\gfr$, and let $\ell$ denote
the Dynkin (multi-)index of the inclusion homomorphism $\pfr \subset \gfr$.
We say that $\pfr$ is a \textit{conformal subalgebra} of $\gfr$ at level
$k$ if $c(\pfr,k\ell) = c(\gfr,k)$.
\end{defi}

The equality $c(\pfr,k\ell)=c(\gfr,k)$ in the above definition in fact holds only if $k= 1$. Classification of
conformal pairs
are given in \cite{BB} and \cite{SW}.  We recall (see \cite{KW}) that, since  $\pfr$ is semisimple,
$\pfr \subset \gfr$ is a conformal subalgebra is equivalent to the statement that
any irreducible $\hat{\gfr}$-module $\cH_\lambda(\gfr)$ of level $1$ decomposes
into a finite  sum of  irreducible $\hat{\pfr}$-modules of level $\ell$.
\bigskip

A fundamental property of a conformal
embedding  $\pfr \subset \gfr$ is the following.

\begin{prop}[\cite{KW} Proposition 3.2 (c)] \label{cosetmodule}
If $\pfr \subset \gfr$ is a conformal embedding, then
$$L^\gfr_n = L^\pfr_n \in \End(\cH_{\hat{\lambda}}(\gfr)) \qquad \text{for all}  \ n \in \ZZ. $$
\end{prop}

In the above statement, we warn that (though implicit in the notation), the operator $L^\gfr_n$ is acting
at level $k=1$ and $L^\pfr_n$ at level $\ell$.  In all the cases we will consider the 
multi-integer $\ell$ is equal to $1$ (see Table (11)).
\bigskip

\subsection{The pairing on $\cH_\lambda \times \cH_{\lambda^\dagger}$}

First we recall the following lemma

\begin{lem}[\cite{U} Lemma 2.2.12] \label{bilinearform}
There exists a bilinear pairing
$$ ( \cdot | \cdot )_\lambda \ : \ \cH_\lambda \times \cH_{\lambda^{\dagger}} \longrightarrow 
\CC,$$
unique up to a multiplicative constant such that
\begin{equation} \label{eqbl}
(X(n) u | v)_\lambda + (u | X(-n) v)_\lambda = 0,
\end{equation}
for any $X \in \gfr$, $ n \in \ZZ$, $u \in \cH_\lambda$ and $v \in \cH_{\lambda^{\dagger}}$. Moreover the
restriction of this pairing to $\cH_\lambda(m) \times \cH_{\lambda^{\dagger}}(m')$ is zero if $m \not= m'$, and non-degenerate
if $m = m'$.
\end{lem}

Consider the restriction of the pairing $( \cdot | \cdot )_\lambda$ to $\cH_\lambda(0) \times \cH_{\lambda^{\dagger}}(0) =
V_\lambda \times V_{\lambda^{\dagger}}$. By definition $V_{\lambda^{\dagger}} = V_\lambda^*$, so that
$V_\lambda \otimes V_{\lambda^{\dagger}} = \End ( V_\lambda)$. The pairing on $V_\lambda \times V_\lambda^*$
is given by a multiple of the natural evaluation map.  More generally, the pairing $( \cdot | \cdot )_\lambda$
induces for any integer $m$ an isomorphism $\phi_m^{\lambda} : \cH_{\lambda^\dagger}(m)
\stackrel{\sim}{\rightarrow} \cH_\lambda(m)^*$
and therefore a distinguished element, which we denote by 
$\gamma_\lambda(m) \in \cH_{\lambda}(m) \otimes \cH_{\lambda^{\dagger}}(m)$ and which maps to
the identity element in $\End(\cH_\lambda(m))$ under the isomorphism
\begin{equation} \label{sewingelements}
id \otimes \phi_m^{\lambda}: \cH_{\lambda}(m) \otimes \cH_{\lambda^{\dagger}}(m) \stackrel{\sim}{\longrightarrow} \End(\cH_\lambda(m)),
\qquad \gamma_\lambda(m) \mapsto \mathrm{Id}_{\cH_\lambda(m)}.
\end{equation}

Note that the family $\{ \gamma_\lambda(m) \}_{m \in \ZZ_+}$ is uniquely defined up to a multiplicative constant. More precisely, if we multiply $( \cdot | \cdot )_\lambda$
with $\alpha \in \CC^*$ the family $\{ \gamma_\lambda(m) \}_{m \in \ZZ_+}$ is transformed
into the family $\{ \frac{1}{\alpha} \gamma_\lambda(m) \}_{m \in \ZZ_+}$.

\bigskip

Consider a conformal embedding $\pfr \subset \gfr$ and an integrable $\hat{\gfr}$-module $\cH_\lambda(\gfr)$ of
level one. Then we have a decomposition as $\hat{\pfr}$-module
$$ \cH_\lambda(\gfr) = \bigoplus_{\hat{\mu} \in B(\lambda)} M(\hat{\mu}, \lambda) \otimes \cH_{\hat{\mu}}(\pfr),$$
where $B(\lambda)$ is a finite subset of $\widehat{P_{\ell}(\pfr)}$  (see \cite{K} equation
(13.14.6)) and the $M(\hat{\mu},\lambda)$ are
finite-dimensional vector spaces (\cite{KW} section 1.6). The integer $\dim M(\hat{\mu}, \lambda)$ is the multiplicity of the representation
$\cH_{\hat{\mu}}(\pfr)$ in $\cH_\lambda(\gfr)$. Note that the weights $\hat{\mu} \in B(\lambda)$ do not necessarily lie in
$P_{\ell}(\pfr)$. We can write $\hat{\mu} =  \mu + \sum_i \ell_i \Lambda_0^{(i)} - n_{\mu} \delta$.

\bigskip

Using Proposition \ref{cosetmodule} we deduce an equality between the trace anomalies
\begin{equation} \label{tracepg}
\Delta_\lambda(\gfr) = \Delta_{\hat{\mu}}(\pfr) \qquad \text{for any}  \ \hat{\mu} \in B(\lambda).
\end{equation}
Moreover by \eqref{traceanomaly} we also have
$\Delta_{\hat{\mu}}(\pfr) = \Delta_\mu(\pfr) - n_{\mu}$. Since $-n_{\mu}$ is the $d$-eigenvalue of the highest weight vector $v_\mu \in \cH_{\hat{\mu}}(\pfr) \subset \cH_\lambda(\gfr)$ and since all $d$-eigenvalues of $\cH_\lambda(\gfr)$ are
negative integers, we conclude that $n_{\mu} \in \ZZ_+$.

\begin{rem} \label{nmuzero}
{\em We immediately deduce from the above that $n_\mu = 0$ if and only if the $\pfr$-module
$V_\mu = \cH_\mu(\pfr) (0)$ appears in the decomposition into irreducible $\pfr$-modules of the $\gfr$-module $V_\lambda = \cH_\lambda(\gfr) (0)$.}
\end{rem}

Thus we conclude that given $\mu \in P_{\ell}(\pfr)$ there exists at most one $\hat{\mu} \in B(\lambda)$ --- since $n_\mu$ is
given by the difference of the trace anomalies. So we will write
$$\mult_\lambda(\mu, \pfr) := \dim M(\hat{\mu}, \lambda)$$
for the multiplicity of occurence of $\cH_\mu(\pfr)$ in $\cH_\lambda(\gfr)$ as $[\hat{\pfr},\hat{\pfr}]$-modules.

\begin{prop} \label{weightsdual}
We have the equality
$$B(\lambda^{\dagger}) = B(\lambda)^{\dagger} := \{ \hat{\mu} \in \widehat{P_{\ell}(\pfr)} \ | \ \hat{\mu}^{\dagger} \in B(\lambda) \}.$$
Moreover $\mult_\lambda(\mu, \pfr) = \mult_{\lambda^\dagger}(\mu^\dagger, \pfr)$.
\end{prop}

\begin{proof}
For $\hat{\nu} \in \widehat{P_k(\gfr)}$ and $\lambda \in P_k(\gfr)$ we will denote by $V^\lambda_{\hat{\nu}} \subset
\cH_\lambda(\gfr)$ the weight space of the $\hat{\gfr}$-module $\cH_\lambda(\gfr)$ associated to the weight $\hat{\nu}$.
It follows from relation \eqref{L0} that $V^\lambda_{\hat{\nu}} \subset \cH_\lambda(\gfr)(m)$ if and only if
$\hat{\nu}$ is of the form $\hat{\nu} = \nu + k\Lambda_0 - m \delta$. By Lemma \ref{bilinearform} we know that
$\cH_\lambda(\gfr)(m)$ and $\cH_{\lambda^\dagger}(\gfr)(m)$ are dual spaces, and it follows from
relation \eqref{eqbl} that the weight spaces
$$ V^\lambda_{\nu + k\Lambda_0 - m \delta} \subset \cH_\lambda(\gfr)(m) \qquad \text{and} \qquad
V^{\lambda^\dagger}_{-\nu + k\Lambda_0 - m \delta} \subset \cH_{\lambda^\dagger}(\gfr)(m)$$
are dual to each other. Hence their dimensions coincide
$$ \mult_{\cH_\lambda}(\nu + k\Lambda_0 - m \delta) = \mult_{\cH_{\lambda^\dagger}}(-\nu + k\Lambda_0 - m \delta). $$

Consider $\hat{\mu} \in B(\lambda)$. Then by \eqref{tracepg} we have
$\Delta_\lambda(\gfr) = \Delta_{\hat{\mu}}(\pfr)$ and since the bilinear form $(\cdot , \cdot )$ is
invariant under the finite Weyl group  and $w_0(\rho) = - \rho$ we also have
$\Delta_{\lambda^\dagger}(\gfr) = \Delta_{\hat{\mu}^\dagger}(\pfr)$. By \cite{K} Proposition 12.11, there
exists a weight $\hat{\nu}$ of $\cH_\lambda(\gfr)$ such that $\hat{\nu}_{|\hat{\hfr}_\pfr} = \hat{\mu}$.
By Proposition \ref{Weylgroup} there exists an element $\tilde{w} \in W(\gfr)$ which restricts
to $w_0^\pfr \in W(\pfr)$ and by \cite{K} Proposition 10.1,
$$ \mult_{\cH_{\lambda^\dagger}}(-\nu + k\Lambda_0 - m \delta) = \mult_{\cH_{\lambda^\dagger}}
(-\tilde{w}(\nu) + k\Lambda_0 - m \delta) .$$
But $-\tilde{w}(\nu) + k\Lambda_0 - m \delta_{|\hat{\hfr}_\pfr} = \hat{\mu}^\dagger$, therefore, using
\cite{K} Proposition 12.11 once more, we obtain that $\hat{\mu}^\dagger \in B(\lambda^\dagger)$. The same reasoning combined with
\cite{K} formula (12.11.1) shows that $\mult_\lambda(\mu, \pfr) = \mult_{\lambda^\dagger}(\mu^\dagger, \pfr)$.
\end{proof}

\subsection{Conformal pairs and the sewing procedure}

Let $q$ be a formal variable. Given $\lambda \in P_k(\gfr)$ we define the element (cf. section 3.3)
$$\tilde{\gamma}_\lambda := \sum_{m=0}^\infty \gamma_\lambda(m) q^m \in \cH_\lambda(\gfr)\otimes \cH_{\lambda^\dagger}(\gfr) [[q]].$$

Note that $\tilde{\gamma}_\lambda$ is well-defined up to a multiplicative constant, and that
the choice of a bilinear pairing $( \cdot | \cdot )_\lambda$ on $\cH_\lambda \times 
\cH_{\lambda^\dagger}$ introduced in Lemma \ref{bilinearform} uniquely determines  
$\tilde{\gamma}_\lambda$.

\bigskip

Let $\pfr \subset \gfr$ be a conformal subalgebra. It follows from Proposition  
\ref{weightsdual} that the decompositions of $\cH_\lambda(\gfr)$ and
$\cH_{\lambda^\dagger}(\gfr)$ into irreducible $\hat{\pfr}$-modules are of the
form
\begin{equation} \label{decHLHLDUAL}
\cH_\lambda(\gfr) = \bigoplus_{\hat{\mu} \in B(\lambda)} \cH_{\hat{\mu}}(\pfr),
\qquad
\cH_{\lambda^\dagger}(\gfr) = \bigoplus_{\hat{\mu} \in B(\lambda)} \cH_{\hat{\mu}^\dagger}(\pfr).
\end{equation}
Here we suppose for simplicity that $\mult_\lambda(\mu, \pfr) = 1$ for all
$\hat{\mu} \in B(\lambda)$, which will be the case of all our examples (see Table (\ref{tabledecomposition}) below). We start with decomposing the bilinear form $(\cdot | \cdot 
)_\lambda$ introduced in Lemma \ref{bilinearform} with respect to the direct sums
\eqref{decHLHLDUAL}.

\begin{lem}
Given $(\hat{\mu}, \hat{\nu}) \in B(\lambda) \times B(\lambda^\dagger)$ the 
restriction of the bilinear form $(\cdot | \cdot )_\lambda$ to
the direct summand $\cH_{\hat{\mu}}(\pfr) \times \cH_{\hat{\nu}}(\pfr)$ is
\begin{itemize}
\item zero, if $\hat{\nu} \not= \hat{\mu}^\dagger$.
\item a non-zero multiple of the bilinear form $(\cdot | \cdot)_\mu$ on
$\cH_{\mu}(\pfr) \times \cH_{\mu^\dagger}(\pfr)$, if $\hat{\nu} = \hat{\mu}^\dagger$.
\end{itemize}
\end{lem}

\begin{proof}
We consider the decomposition into $d$-eigenspaces
$$ \cH_\lambda(\gfr)(m) = \bigoplus_{\hat{\mu} \in B(\lambda)} \cH_{\hat{\mu}}(\pfr) (m) .$$
Note that $\cH_{\hat{\mu}}(\pfr)(m) = \cH_\mu(\pfr)(m - n_\mu)$ for all $m \in \ZZ_+$ and that $\cH_\mu(\pfr)(l) = \{ 0 \}$ for $l < 0$. First we observe that the restriction of 
$(\cdot | \cdot )_\lambda$ to $\cH_{\hat{\mu}}(\pfr) \times \cH_{\hat{\nu}}(\pfr)$
is determined by values on the finite-dimensional subspace
$$ V_\mu \times V_\nu  = \cH_\mu(\pfr)(0) \times \cH_\nu (\pfr)(0) = 
 \cH_{\hat{\mu}}(\pfr)(n_\mu) \times \cH_{\hat{\nu}} (\pfr)(n_\nu).$$
This follows from the fact that $(\cdot | \cdot )_\lambda$ satisfies relation \eqref{eqbl}
for $X \in \pfr, u \in \cH_{\hat{\mu}}(\pfr)$ and $v \in \cH_{\hat{\nu}}(\pfr)$, which enables
one to reconstruct the bilinear form on $\cH_{\hat{\mu}}(\pfr) \times \cH_{\hat{\nu}}(\pfr)$
from its values on $V_\mu \times V_\nu$ --- see e.g. \cite{U} Lemma 2.2.12.

\bigskip

The restriction of $( \cdot | \cdot )_\lambda$ to $V_\mu \times V_\nu$ induces a 
linear map $\varphi : V_\mu \rightarrow (V_\nu)^*$ and relation \eqref{eqbl} written
for $X \in \pfr, n= 0, u \in V_\mu$ and $v \in V_\nu$ implies that $\varphi$ is a 
$\pfr$-module homomorphism (here we endow the vector space $(V_\nu)^*$ with the 
structure of a left $\pfr$-module). Hence in the case $\hat{\nu} \not= \hat{\mu}^\dagger$,
we deduce that $\varphi = 0$ and therefore that the restriction of $( \cdot | \cdot )_\lambda$
to $\cH_{\hat{\mu}}(\pfr) \times \cH_{\hat{\nu}}(\pfr)$ is zero. In the case 
$\hat{\nu} = \hat{\mu}^\dagger$, we obtain by Schur's lemma that $\varphi : V_\mu 
\rightarrow (V_{\mu^\dagger})^* = V_\mu$ is a homothety and that the restriction of
$( \cdot | \cdot )_\lambda$ to $\cH_{\hat{\mu}}(\pfr) \times \cH_{\hat{\mu}^\dagger}(\pfr)$
is a multiple of $(\cdot | \cdot)_\mu$.  It now remains to check that this multiple is non-zero. In order to show that, we consider
for $m$ large the isomorphism $\phi_m^{\lambda}$ induced by $(\cdot | \cdot )_\lambda$
$$\cH_{\lambda^\dagger}(m) = \bigoplus_{\hat{\mu} \in B(\lambda)} \cH_{\hat{\mu}^\dagger}(m)
\longrightarrow
\cH_\lambda(m)^* = \bigoplus_{\hat{\mu} \in B(\lambda)} \cH_{\hat{\mu}}(m)^*.$$
By the preceding arguments, $\phi_m^{\lambda}$ is block-diagonal and decomposes as 
$\phi_m^{\lambda}=\sum_{\hat{\mu} \in B(\lambda)}\phi_m^{\mu}$.
Hence all restrictions 
$\phi_m^{\mu}:  \cH_{\hat{\mu}^\dagger}(m)  \rightarrow
\cH_{\hat{\mu}}(m)^*$  for $\hat{\mu} \in B(\lambda)$ are isomorphisms and we are done.
\end{proof}

\begin{rem} \label{normalizegammas}
{\em If we fix a pairing $(\cdot | \cdot )_\lambda$ on $\cH_\lambda \times \cH_{\lambda^\dagger}$, 
we also obtain by restriction a pairing $(\cdot | \cdot)_\mu$ on 
$\cH_\mu \times \cH_{\mu^\dagger}$ for any $\hat{\mu} \in B(\lambda)$.  Thus we obtain
uniquely defined elements $\tilde{\gamma}_\lambda$ and $\tilde{\gamma}_\mu$  for
$\hat{\mu} \in B(\lambda)$.}
\end{rem}

\begin{prop} \label{decgammalambda}
With the choices made in Remark \ref{normalizegammas} we have a decomposition in $\cH_\lambda(\gfr) \otimes \cH_{\lambda^\dagger}(\gfr) [[q]]$
$$ \tilde{\gamma}_\lambda = \sum_{\hat{\mu} \in B(\lambda)} q^{n_\mu} \tilde{\gamma}_\mu,$$
where the positive integer $n_\mu$ equals $\Delta_\mu(\pfr) - \Delta_\lambda(\gfr)$.
\end{prop}

\begin{proof}
As in the proof of Lemma 3.6, 
we consider the decomposition into $d$-eigenspaces
$$\cH_\lambda(\gfr)(m) = \bigoplus_{\hat{\mu} \in B(\lambda)} \cH_{\hat{\mu}}(\pfr) (m)$$ and 
note that $\cH_{\hat{\mu}}(\pfr)(m) = \cH_\mu(\pfr)(m - n_\mu)$ for all $m \in \ZZ_+$. 
The identity transformation of $\cH_\lambda(\gfr)(m)$ obviously decomposes as
\begin{equation}\label{iddecomp}
\mathrm{Id}_{\cH_\lambda(\gfr)(m)} = \sum_{\hat{\mu} \in B(\lambda)} \mathrm{Id}_{\cH_{\mu}(\pfr)(m - n_\mu)} \qquad
\end{equation}
for any $\ m \in \ZZ_+$.

Since by Lemma 3.6 the pairing $( \cdot | \cdot )_\lambda$ on $\cH_\lambda \times
\cH_{\lambda^\dagger}$  restricts to the pairing $( \cdot | \cdot )_\mu$ on $\cH_\mu \times \cH_{\mu^\dagger}$
for any $\hat{\mu}\in B(\lambda)$, we obtain after applying the inverse of the isomorphism 
$id \otimes \phi_m^{\lambda}=\sum_{\hat{\mu} \in B(\lambda)}id \otimes \phi_m^{\mu}$ (see (\ref{sewingelements})) 
to the equality (\ref{iddecomp}),
$$ \gamma_\lambda(m) = \sum_{\hat{\mu} \in B(\lambda)} \gamma_\mu(m - n_\mu)  \qquad
\text{for any} \ m \in \ZZ_+.$$
Here we put $\gamma_\mu(l) = 0$ if $l <0$.
Multiplying with $q^m$ and summing over $\ZZ_+$ gives the relation
\begin{eqnarray*}
\tilde{\gamma}_\lambda = \sum_{m = 0}^\infty \gamma_\lambda(m) q^m  & = & \sum_{\hat{\mu} \in B(\lambda)}
\sum_{m \geq n_\mu} \gamma_\mu(m - n_\mu) q^{m - n_\mu} q^{n_\mu} \\
& =  & \sum_{\hat{\mu} \in B(\lambda)} q^{n_\mu} \tilde{\gamma}_\mu.
\end{eqnarray*}

\end{proof}

The following table is extracted from \cite{KS} page 2235 and gives for any conformal subalgebra $\pfr$ of
Table \eqref{confE8} the decomposition of the basic
representation $\cH_0(\mathfrak{e}_8)$ as $[\hat{\pfr}, \hat{\pfr}]$-module, in particular its
Dynkin (multi-)index $\ell$, its subset $B(0)$ and the action of the
involution $\mu \mapsto \mu^\dagger$. We use Bourbaki's notation for the fundamental weights of a simple
Lie algebra.

\begin{equation} \label{tabledecomposition}
\begin{tabular}{|c||c|c|c|}
        \hline
             $\pfr$ & $\ell$ & B(0)  & $\varpi \mapsto \varpi^\dagger$  \\
            \hline
            \hline
            $\mathfrak{so}(16)$ & 1 &  $\{0, \varpi_7 \}$ & $\varpi_7^\dagger = \varpi_7$ \\
     \hline
      $\mathfrak{sl}(9)$ & 1 & $\{ 0,\varpi_3, \varpi_6 \}$ & $\varpi_3^\dagger = \varpi_6$ \\
      \hline
      $\mathfrak{sl}(5) \oplus \mathfrak{sl}(5)$ & (1,1) & $\{ (0,0) ,(\varpi_1, \varpi_2),
      (\varpi_2, \varpi_4),(\varpi_3, \varpi_1),(\varpi_4, \varpi_3) \}$ &
      $(\varpi_i, \varpi_j)^\dagger = (\varpi_{5-i}, \varpi_{5-j})$ \\
      \hline
       $\mathfrak{sl}(3) \oplus \mathfrak{e}_6$ & (1,1) & $\{ (0,0) ,(\varpi_1, \varpi_1),
      (\varpi_2, \varpi_6) \}$ &
      $(\varpi_1, \varpi_1)^\dagger = (\varpi_2, \varpi_6)$ \\
      \hline
      $\mathfrak{sl}(2) \oplus \mathfrak{e}_7 $ & (1,1) & $\{ (0,0) ,(\varpi_1, \varpi_7) \}$ &
      $(\varpi_1, \varpi_7)^\dagger = (\varpi_1, \varpi_7)$ \\
      \hline
      $\mathfrak{g}_2 \oplus \mathfrak{f}_4$ & (1,1) & $\{ (0,0) ,(\varpi_1, \varpi_4) \}$ &
      $(\varpi_1, \varpi_4)^\dagger = (\varpi_1, \varpi_4) $ \\
      \hline

\end{tabular}
\end{equation}

We recall that if $\pfr = \mathfrak{a} \oplus \mathfrak{b}$, then we have the
decomposition as $[\hat{\pfr}, \hat{\pfr}]$-module
$$\cH_0(\mathfrak{e}_8) = \bigoplus_{(\lambda_1, \lambda_2) \in B(0)} \cH_{\lambda_1}(\mathfrak{a})
\otimes \cH_{\lambda_2}(\mathfrak{b}).$$
We also note that, for any $\pfr$ in Table (\ref{tabledecomposition}), the trace
anomalies $\Delta_\mu(\pfr) = n_\mu$ are equal to $1$ for all non-zero $\mu \in B(0)$.

\section{Conformal pairs: factorization rules}

In this section we describe how the factorization rules behave under conformal embeddings. First we recall the factorization rules. Let $X_0$ be a nodal (not necessarily irreducible) curve with one node $x_0$. We call $\tilde{X}$ the normalization
of $X_0$ with $\pi: \tilde{X} \rightarrow X_0$ and $\pi^{-1}(x_0) = \{ a,b \} $.

\begin{prop}[Factorization rules, \cite{U} Theorem 4.4.9] \label{fact}
Let $\gfr$ be a semi-simple Lie algebra and $X_0$ a nodal curve with $s$ marked points with labels
$\vec{\lambda} \in P_k(\gfr)^s$.
There is an isomorphism
$$ \bigoplus_{\lambda \in P_k(\gfr)} \cV^\dagger_{\vec{\lambda}, \lambda, \lambda^\dagger} (\tilde{X}, \gfr)
\stackrel{\oplus \iota_\lambda}{\longrightarrow} \cV^\dagger_{\vec{\lambda}}(X_0,\gfr).$$
\end{prop}

\begin{rem} \label{normalizeiota}
{\em For any $\lambda \in P_k(\gfr)$ the linear map $\iota_\lambda$ is only defined up to
homothety. More precisely, $\iota_\lambda$ depends on the choice of an element $\tilde{\gamma}_\lambda$ (or,
equivalently, of its degree zero piece $\gamma_\lambda(0)$), which are only defined up to
a multiplicative constant --- see equation \eqref{defiotalambda} below.}
\end{rem}

\bigskip

We will denote by $\cO$ the ring of formal power series $\CC[[q]]$ and by $K = \CC((q))$ its field of
fraction. We consider a family of curves $\cX$ over $\Spec \ \cO$ such that its special fiber $\cX_0$ is a nodal
curve $X_0$ over $\CC$ and its generic fiber $\cX_K$ a smooth curve over the field  $K$. Consider the
sheaf of conformal blocks $\cV^\dagger_{\vec{\lambda}}(\cX, \gfr)$ for the family $\cX$, which is
an $\cO$-module of finite type. Moreover since the formation of the sheaf of conformal blocks commutes
with base change, the fibre $\cV^\dagger_{\vec{\lambda}}(\cX, \gfr)_0$ over the closed point $0 \in \Spec \ \cO$ of
$ \cV^\dagger_{\vec{\lambda}}(\cX, \gfr)$ coincides with $ \cV^\dagger_{\vec{\lambda}}(X_0, \gfr)$. We thus obtain a
restriction map
$$ r_0 : \cV^\dagger_{\vec{\lambda}}(\cX, \gfr) \longrightarrow \cV^\dagger_{\vec{\lambda}}(X_0, \gfr).$$
On the other hand, there exists for every $\lambda \in P_k(\gfr)$ a $\CC$-linear map --- the so-called sewing procedure, see
\cite{U} formula (4.4.3) and Lemma 4.4.5
$$ \mathfrak{s}_\lambda : \cV^\dagger_{\vec{\lambda}, \lambda, \lambda^\dagger}(\tilde{X}, \gfr)  \longrightarrow
 \cV^\dagger_{\vec{\lambda}}(\cX, \gfr), \qquad \psi_\lambda \mapsto \widetilde{\psi_\lambda} :=
 \mathfrak{s}_\lambda(\psi_\lambda).$$
The linear maps $\iota_\lambda$ and $\mathfrak{s}_\lambda$ are defined as follows: for $\psi_\lambda
\in \cV^\dagger_{\vec{\lambda}, \lambda, \lambda^\dagger}(\tilde{X}, \gfr)$
\begin{equation} \label{defiotalambda}
\langle \iota_\lambda (\psi_\lambda) | u \rangle := \langle \psi_\lambda | u \otimes \gamma_\lambda(0) \rangle \in
\CC \qquad
\text{and} \qquad \langle \widetilde{\psi_\lambda} | \tilde{u} \rangle := \langle \psi_\lambda |
\tilde{u} \otimes \tilde{\gamma}_\lambda \rangle \in \cO
\end{equation}
for any vectors $u \in \cH_{\vec{\lambda}}$ and $\tilde{u} \in \cH_{\vec{\lambda}}[[q]]$. We recall (see \cite{U}
Lemma 4.4.6) that $\cV^\dagger_{\vec{\lambda}}(\cX, \gfr)$ identifies with the subset of linear forms
in $\cH^\dagger_{\vec{\lambda}}[[q]]$ satisfying the formal gauge condition. It is clear from these definitions that
the map $\mathfrak{s}_\lambda$ is an extension of $\iota_\lambda$, i.e.,
$\iota_\lambda = r_0 \circ \mathfrak{s}_\lambda$.

\bigskip

We consider now a conformal embedding $\pfr \subset \gfr$. We assume that all level one representations
$\cH_\lambda(\gfr)$ decompose with multiplicities one, i.e. $\mult_\lambda(\mu, \pfr) = 1$ for all
$\hat{\mu} \in B(\lambda)$. For $\vec{\lambda} = (\lambda_1, \ldots, \lambda_s)$ define
$B(\vec{\lambda}) = B(\lambda_1) \times \cdots \times B(\lambda_s)$. Consider $\vec{\hat{\mu}} \in B(\vec{\lambda})$ and the
corresponding inclusion
$$ \cH_{\vec{\hat{\mu}}}(\pfr) \hookrightarrow \cH_{\vec{\lambda}}(\gfr).$$
The gauge condition is preserved under restriction of linear forms to
$\cH_{\vec{\hat{\mu}}}(\pfr)$ --- see \cite{NT} formula (2.9), so that we obtain an $\cO$-linear map
\begin{equation} \label{mapalpha}
\alpha:  \cV^\dagger_{\vec{\lambda}}(\cX, \gfr) \longrightarrow \cV^\dagger_{\vec{\mu}}(\cX, \pfr).
\end{equation}
The restriction of $\alpha$ to $0 \in \Spec \ \cO$ gives  a $\CC$-linear map
$$\alpha_0:  \cV^\dagger_{\vec{\lambda}}(X_0, \gfr) \longrightarrow \cV^\dagger_{\vec{\mu}}(X_0, \pfr). $$

The next result will describe how $\alpha_0$ decomposes in the direct sums given by the 
factorization rules (Proposition \ref{fact}) on both sides. More precisely, we consider the
composite map

$$\tilde{\alpha} : \bigoplus_{\lambda \in P_1(\gfr)} \cV^\dagger_{\vec{\lambda}, \lambda, \lambda^\dagger} (\tilde{X}, \gfr)
\stackrel{\oplus \iota_\lambda}{\longrightarrow} \cV^\dagger_{\vec{\lambda}}(X_0, \gfr)  \stackrel{\alpha_0}{\longrightarrow} \cV^\dagger_{\vec{\mu}}(X_0, \pfr) \stackrel{(\oplus \iota_\mu)^{-1}}{\longrightarrow} 
\bigoplus_{\mu \in P_{\ell}(\pfr)} \cV^\dagger_{\vec{\mu}, \mu, \mu^\dagger} (\tilde{X}, \pfr) $$
and define for any pair $(\lambda, \mu) \in P_1(\gfr) \times P_\ell(\pfr)$ the linear map
$\alpha_0^{\lambda,\mu} : \cV^\dagger_{\vec{\lambda}, \lambda, \lambda^\dagger}(\tilde{X}, \gfr) \rightarrow
\cV^\dagger_{\vec{\mu}, \mu, \mu^\dagger} (\tilde{X}, \pfr)$ to be the $(\lambda,\mu)$-component of $\tilde{\alpha}$.
In particular, we have for any $\lambda \in P_1(\gfr)$, $\sum_{\mu \in P_{\ell}(\pfr)} \iota_\mu \circ \alpha_0^{\lambda,\mu} =
\alpha_0 \circ \iota_\lambda$.

\begin{rem}
{\em Since $\iota_\lambda$ and $\iota_\mu$ are only defined up to homothety, the linear maps $\alpha_0^{\lambda,\mu}$ are only
defined up to homothety.}
\end{rem}

\begin{prop} \label{mainresult}
Given a pair $(\lambda, \mu) \in P_1(\gfr) \times P_\ell(\pfr)$, the linear map $\alpha_0^{\lambda,\mu}$ 
\begin{itemize}
\item is identically zero, if $\hat{\mu} \notin B(\lambda)$ or $\Delta_\mu (\pfr) \not= \Delta_\lambda(\gfr)$.
\item is, up to non-zero homothety, induced by the natural inclusion $\cH_{\vec{\hat{\mu}}}(\pfr) \otimes \cH_{\hat{\mu}} (\pfr) \otimes
\cH_{\hat{\mu}^\dagger}(\pfr) \hookrightarrow \cH_{\vec{\lambda}}(\gfr) \otimes  \cH_{\lambda}(\gfr)
\otimes \cH_{\lambda^\dagger}(\gfr)$, if $\hat{\mu} \in B(\lambda)$ and  $\Delta_\mu (\pfr) = \Delta_\lambda(\gfr)$.
\end{itemize}
\end{prop}

\begin{proof}
We fix a weight $\lambda \in P_1(\gfr)$. In order to lift the $\CC^*$-indeterminacy in the definition of $\alpha_0^{\lambda, \mu}$, we 
fix a bilinear pairing $(\cdot | \cdot )_\lambda$ on $\cH_\lambda \times \cH_{\lambda^\dagger}$, which by Remarks
\ref{normalizegammas} and \ref{normalizeiota} determines the elements $\tilde{\gamma}_\lambda$ and $\tilde{\gamma}_\mu$,
hence the linear maps $\iota_\lambda$ and $\iota_\mu$ for any $\hat{\mu} \in B(\lambda)$; we choose any $\iota_\mu$ for $\hat{\mu}
\notin B(\lambda)$.

\bigskip

In order to compute the decomposition of $\alpha_0 (\iota_\lambda (\psi_\lambda))$ for an element
$\psi_\lambda \in \cV^\dagger_{\vec{\lambda}, \lambda, \lambda^\dagger} (\tilde{X}, \gfr)$ in the direct sum $\bigoplus_{\mu \in P_{\ell}(\pfr)} \cV^\dagger_{\vec{\mu}, \mu, \mu^\dagger} (\tilde{X}, \pfr)$, we will first decompose  the extension $\alpha (\widetilde{\psi_\lambda})$ and then
restrict to the special fiber. In fact, we have $ \alpha_0 (\iota_\lambda (\psi_\lambda)) = r_0(\alpha (\widetilde{\psi_\lambda}))$.
On the other hand, using Proposition \ref{decgammalambda} and the definition \eqref{defiotalambda} of $\widetilde{\psi_\lambda}$, we easily
obtain the following equalities in $\cO$, which hold for any $\tilde{u} \in \cH_{\vec{\hat{\mu}}}(\pfr)[[q]] \hookrightarrow 
\cH_{\vec{\lambda}}(\gfr)[[q]]$
\begin{equation} \label{eqpsilambda}
\langle \alpha(\widetilde{\psi_\lambda}) | \tilde{u} \rangle = \langle \widetilde{\psi_\lambda} | \tilde{u} \rangle
 =  \langle \psi_\lambda | \tilde{u} \otimes \tilde{\gamma}_\lambda \rangle 
=  \sum_{\hat{\mu} \in B(\lambda)} q^{n_\mu} \langle \psi_\lambda | \tilde{u} \otimes \tilde{\gamma}_\mu \rangle.
\end{equation}
Now we restrict this expression to the special fiber, i.e., we put $q = 0$: we write $\tilde{u} = \sum_{m \geq 0} u(m)q^m$
and $\tilde{\gamma}_\mu = \sum_{m \geq 0} \gamma_\mu(m) q^m$, so that the evaluation at $q=0$ of the first term of 
\eqref{eqpsilambda} gives
\begin{equation}\label{LHS}
\langle \alpha(\widetilde{\psi_\lambda}) | \tilde{u} \rangle (0) = \langle r_0 (\alpha (\widetilde{\psi_\lambda})) | u(0)
\rangle = \langle \alpha_0 (\iota_\lambda( \psi_\lambda )) | u(0) \rangle.
\end{equation}
On the other hand the evaluation at $q=0$ of the last term of \eqref{eqpsilambda} annihilates all terms corresponding
to $\hat{\mu} \in B(\lambda)$ with $n_\mu \not= 0$ and, since $[\tilde{u} \otimes \tilde{\gamma}_\mu](0) = u(0)
\otimes \gamma_\mu(0)$, the last term becomes
\begin{equation} \label{RHS}
\sum_{\hat{\mu} \in B(\lambda) \atop{n_\mu = 0}} \langle \psi_\lambda | u(0) \otimes \gamma_\mu (0) \rangle =
\sum_{\hat{\mu} \in B(\lambda) \atop{n_\mu = 0}} \langle \mathrm{res}_\mu(\psi_\lambda) | u(0) \otimes \gamma_\mu (0) \rangle =
\sum_{\hat{\mu} \in B(\lambda) \atop{n_\mu = 0}} \langle \iota_\mu (\mathrm{res}_\mu(\psi_\lambda)) | u(0) \rangle, 
\end{equation}
where $\mathrm{res}_\mu(\psi_\lambda)$ denotes the restriction of the linear form $\psi_\lambda$ to the
subspace $\cH_{\vec{\hat{\mu}}}(\pfr) \otimes \cH_{\hat{\mu}} (\pfr) \otimes
\cH_{\hat{\mu}^\dagger}(\pfr)$. Since the equality between \eqref{LHS} and \eqref{RHS} holds for any vector
$u(0) \in \cH_{\vec{\hat{\mu}}}(\pfr)$, we obtain the following equality in $\cV^\dagger_{\vec{\mu}}(X_0, \pfr)$
$$ \alpha_0 (\iota_\lambda( \psi_\lambda )) = \sum_{\hat{\mu} \in B(\lambda) \atop{n_\mu = 0}} \iota_\mu( \mathrm{res}_\mu (\psi_\lambda) ).$$
Projecting onto the subspace $\mathrm{im}(\iota_\mu) \subset \cV^\dagger_{\vec{\mu}}(X_0, \pfr)$ leads to
$\alpha_0^{\lambda,\mu} = 0$ if $\hat{\mu} \notin B(\lambda)$ or $n_\mu \not= 0$, and $\alpha_0^{\lambda,\mu} (\psi_\lambda) = 
\mathrm{res}_\mu (\psi_\lambda)$ if $\hat{\mu} \in B(\lambda)$ and $n_\mu = 0$.
\end{proof}

\section{Proof of  Theorem 1}

\subsection{Proof of part (i)}
First of all, we will restate Theorem 1 in terms of spaces of conformal blocks. By \cite{B} Proposition 5.2 the linear
map $\phi_P$ identifies under the ``Verlinde" isomorphism \eqref{verlindeiso} with the linear map
$$ \phi_{X, \mathfrak{p}} :  \cV^\dagger_0(X, \mathfrak{e}_8) \longrightarrow \cV^\dagger_0(X, \pfr),$$
induced by the inclusion of the basic $\hat{\pfr}$-module $\cH_0(\pfr)$ into the basic $\hat{\mathfrak{e}_8}$-module  $\cH_0(\mathfrak{e}_8)$. We recall that we choose a point
$p \in X$ and we denote $U = X \setminus \{ p \}$. It is therefore equivalent to show that
$\phi_{X, \mathfrak{p}}$
is non-zero for any smooth curve $X$.

\bigskip

We then observe that the rank of the linear map $\phi_{X, \mathfrak{p}}$ is constant when the smooth curve
$X$ varies by \cite{B} Proposition 5.8 and Lemma A.1. It is therefore sufficient to show that there exists
a smooth curve $X$ for which the map $\phi_{X, \mathfrak{p}}$ is non-zero. We will prove that by induction on the
genus $g$ of $X$.

\bigskip

The case $g=0$ is easily seen as follows. Over the projective line we have $\mathfrak{g}(U) \cdot \cH_0(\gfr) =
\oplus_{m > 0} \cH_0(\gfr)(m)$ for any semi-simple Lie algebra $\gfr$. Hence the one-dimensional space
of covacua $\cV_0 (\PP^1, \gfr)$ is generated by the image under the projection
$$ \cH_0(\gfr) \longrightarrow \cV_0(\PP^1, \gfr) = \cH_0(\gfr) / \mathfrak{g}(U) \cdot \cH_0(\gfr) $$
of the highest weight vector  $v_0(\gfr) \in \cH_0(\gfr)(0) = V_0 = \CC$. Since the trivial
$\mathfrak{e}_8$-module $V_0 = \CC$ restricts to the trivial $\mathfrak{p}$-module, the highest weight
vector $v_0(\pfr) \in \cH_0(\pfr)(0) = V_0$ coincides with the highest weight vector $v_0(\gfr)$
under the inclusion $\cH_0(\pfr) \hookrightarrow \cH_0(\gfr)$, which implies that $\phi_{\PP^1, \mathfrak{p}}$
is non-zero.

\bigskip

Next, we consider as in section 4 a family $\cX$ of genus $g$ curves parametrized by $\Spec \  \cO$ such
that $\cX_0 = X_0$ is a nodal curve and $\cX_K$ a smooth curve defined over $K = \CC((q))$. We also
consider \eqref{mapalpha} the $\cO$-linear map
$\alpha$ associated to the conformal embedding $\pfr \subset \mathfrak{e}_8$ and the trivial weights $\lambda = 0$ and $\mu = 0$
$$ \alpha : \cV^\dagger_0(\cX, \mathfrak{e}_8) \longrightarrow \cV^\dagger_0(\cX,\pfr).$$
By Proposition \ref{mainresult} and by the induction hypothesis, the restriction $\alpha_0$ of
the map $\alpha$ to the special fiber is non-zero:  in fact, the genus of the normalization
$\tilde{X}$ equals $g-1$ and $\alpha_0$ decomposes, up to non-zero homothety in each direct summand, as follows
$$ \alpha_0 = (\alpha_0^{0,\mu}) :   \cV^\dagger_0(\tilde{X}, \mathfrak{e}_8)
\stackrel{\sim}{\leftarrow} \cV^\dagger_{0,0,0}(\tilde{X}, \mathfrak{e}_8) \lra
\bigoplus_{\mu \in P_{\ell}(\pfr)} \cV^\dagger_{0, \mu, \mu^\dagger}(\tilde{X}, \pfr),$$
where the first isomorphism is the so-called {\em propagation of vacua} isomorphism (see e.g.
\cite{U} Theorem 3.3.1). We note that $P_1(\mathfrak{e}_8) = \{ 0 \}$. By Remark \ref{nmuzero} and
Proposition \ref{mainresult} we have
$\alpha_0^{0, \mu} = 0$ if $\mu \not= 0$ and, again due to the
isomorphism $\cV^\dagger_{0, 0, 0}(\tilde{X}, \pfr) \stackrel{\sim}{\rightarrow}
\cV^\dagger_0(\tilde{X}, \pfr)$, the map $\alpha_0^{0,0}$ is identified with the map $\phi_{\tilde{X}, \mathfrak{p}} :
\cV^\dagger_0(\tilde{X}, \mathfrak{e}_8) \longrightarrow \cV^\dagger_0(\tilde{X}, \pfr)$, which
is non-zero by the induction hypothesis. Hence $\alpha_0$ is also non-zero.

\bigskip

By semi-continuity we conclude that the restriction $\alpha_K = \phi_{\cX_K, \mathfrak{p}} :
\cV^\dagger_0(\cX_K, \mathfrak{e}_8) \longrightarrow \cV^\dagger_0(\cX_K,\pfr)$ of $\alpha$ to the generic fiber is non-zero. Hence, again by \cite{B} Proposition 5.8 and Lemma A.1, the $K$-linear map $\phi_{X_K, \mathfrak{p}}$ is non-zero for any genus $g$ curve $X_K$ defined over the field $K = \CC((q))$. Given a genus $g$ curve $X$ defined
over $\CC$, the result then follows from the equality $\phi_{X, \mathfrak{p}} \otimes_{\CC} K = \phi_{X\otimes_\CC K, \mathfrak{p}}$.

\subsection{Irreducible representations of Heisenberg groups} Before pursuing the proof of Theorem 1, we need
to recall some known facts on Heisenberg groups and their irreducible representations. We consider a semi-simple
simply connected group $P$ with center $Z$ of the following table

\begin{equation} \label{listPwithcenter}
\begin{tabular}{|c||c|c|c|c|c|c|c|}
        \hline
            P & $\Spin(4n)$ &  $\SL(n)$ & $E_6$ & $E_7$ & $\SL(n) \times \SL(n)$ & $\SL(3) \times E_6$ &  $\SL(2) \times E_7$ \\
            \hline
     Z & $\ZZ/2 \ZZ \times \ZZ/2 \ZZ$ & $\ZZ/ n\ZZ$  & $\ZZ/ 3\ZZ$ & $\ZZ/ 2\ZZ$ & $\ZZ/n\ZZ \times \ZZ/n\ZZ$ & $\ZZ/3 \ZZ
     \times \ZZ/3 \ZZ$ & $\ZZ/ 2 \ZZ \times \ZZ/ 2 \ZZ$
     \\
     \hline
\end{tabular}
\end{equation}

\bigskip

The finite abelian group $\cM_X(Z)$ of principal $Z$-bundles acts on $\cM_X(P)$ by twisting $P$-bundles with
$Z$-bundles. Note that $|\cM_X(Z)| = |Z|^{2g}$. We denote by $t_\zeta$ the automorphism of $\cM_X(P)$ induced by the twist with $\zeta \in
\cM_X(Z)$. We introduce the Mumford group associated to the line bundle $\cL_P$
$$ \cG(\cL_P) := \{ (\zeta, \psi) \ | \zeta \in \cM_X(Z) \ \text{and} \ \psi: t^*_\zeta \cL_P
\stackrel{\sim}{\rightarrow} \cL_P \}.$$
The Mumford group is a central extension of the group  $\cM_X(Z)$ by $\CC^*$ and acts via $s \mapsto
\psi (t^*_\zeta (s))$ on the space of global sections $H^0(\cM_X(P), \cL_P)$. Note that
the center $\CC^*$ of $\cG(\cL_P)$ acts by scalar multiplication.

\bigskip

We will need the following results.

\begin{lem}[\cite{F1} Lemma 16] \label{idheis} 
The Mumford group $\cG(\cL_P)$ is isomorphic as an extension to a finite Heisenberg group $\cG(\delta)$. The type
$\delta$ depends only on the center $Z$ and on the genus $g$ of $X$.
\end{lem}

We refer, for example, to \cite{M} page 294 for the definition of the Heisenberg group $\cG(\delta)$.  Under the above identification, 
for all groups $P$ of Table (\ref{listPwithcenter}) the space of global
sections $H^0(\cM_X(P), \cL_P)$ is an irreducible representation of $\cG(\cL_P)$.  This is an immediate consequence
of the fact that there exists a
unique irreducible representation of $\cG(\delta)$ of dimension $|Z|^g$ on
which $\CC^*$ acts as scalar multiplication (see e.g. \cite{M} Proposition 3), and the numerical
identity $\dim H^0(\cM_X(P), \cL_P) = |Z|^{g}$ provided by the Verlinde formula (see Table (\ref{tableVerlinde})).

\bigskip

We remark that for all pairs $(P,N)$ of Table (\ref{listP}), $\cM_X(N)$ is a maximal isotropic subgroup of $\cM_X(Z)$.
Here isotropic means with respect to the standard
symplectic form on $\cM_X(Z)$ induced by the commutators in the Mumford group (see e.g. \cite{M}
page 293). The isotropic subgroup $\cM_X(N)$ is maximal since $|N|=|Z|^{1/2}$, hence $|\cM_X(N)|=|\cM_X(Z)|^{1/2}$. Therefore, there
exists a lift $\cM_X(N) \hookrightarrow \cG(\cL_P)$; in other words, there exists an $\cM_X(N)$-linearization of the line bundle $\cL_P$.

\bigskip

Since $H^0(\cM_X(P), \cL_P)$ is the unique irreducible representation of $\cG(\cL_P)$ of level one, i.e., $\CC^* \subset \cG(\cL_P)$ acts as scalar multiplication, we deduce from Mumford's theory of theta groups (see \cite{M} Proposition 3) the following:

\begin{lem}\label{irrspaceglobal}
For any $\cM_X(N)$-linearization of $\cL_P$, the subspace of $\cM_X(N)$-invariant sections of 
$H^0(\cM_X(P), \cL_P)$ is one-dimensional.
\end{lem}

\subsection{Proof of part (ii)}
It is clear that the image of $\phi$ is contained in the $\cM_X(N)$-invariant subspace.  Hence it suffices to show that the $\cM_X(N)$-invariant
subspace of $H^0(\cM_X(P), \cL_P)$ is one-dimensional, which is precisely the result of Lemma \ref{irrspaceglobal}.

\section{Proof of Theorem 2}

First of all, we recall that by Lemma \ref{irrspaceglobal} the subspace of $\cM_X(N)$-invariant
sections is one-dimensional. The main observation is that, in all three cases, the finite group $N \subset P = A \times B$ projects isomorphically to the centers $Z(A)$ and $Z(B)$ of the 
groups $A$ and $B$ respectively. Hence we obtain canonical isomorphisms $Z(B) 
\stackrel{\sim}{\rightarrow} N \stackrel{\sim}{\rightarrow} Z(A)$, which
induce an isomorphism $\iota : \cM_X(Z(B)) \stackrel{\sim}{\rightarrow} \cM_X(N) \stackrel{\sim}{\rightarrow} \cM_X(Z(A))$. Now we observe that the isomorphism $\iota$ lifts
to an isomorphism
$$ \tilde{\iota} : \cG(\cL_B) \longrightarrow \cG(\cL_A), \qquad x \mapsto 
\tilde{\iota}(x) := \Lambda(x) \cdot x^{-1}.$$
Here $\Lambda$ denotes the composite map $\cG(\cL_B) \rightarrow \cM_X(Z(B)) 
\stackrel{\sim}{\rightarrow} \cM_X(N) \hookrightarrow \cG(\cL_P)$, where 
the first arrow is the canonical projection and the last arrow is the chosen
lift of $\cM_X(N)$ to the Mumford group. In order to show that $\tilde{\iota}$ is a
group homomorphism, we need the fact that the two Mumford groups $\cG(\cL_A)$
and $\cG(\cL_B)$ are commuting subgroups of $\cG(\cL_{A \times B})$. Note that the
restriction of $\tilde{\iota}$ to the center $\CC^*$ is $t \mapsto t^{-1}$.

\bigskip

Under the identification of the Mumford groups $\tilde{\iota} : \cG(\cL_B) \stackrel{\sim}{\rightarrow} \cG(\cL_A)$ the vector space $H^0(\cM_X(A) , \cL_A)^*$ can be viewed
as an irreducible representation of $\cG(\cL_B)$ of level one.
Moreover the $\cM_X(N)$-invariance
of the section $\sigma \in H^0(\cM_X(P), \cL_P)$ translates into the equivariance of 
the linear map
$\sigma : H^0(\cM_X(A) , \cL_A)^* \lra H^0(\cM_X(B), \cL_B)$ under the Mumford group 
$\cG(\cL_B)$.
Since both spaces are irreducible representations of $\cG(\cL_B)$ by section 5.2, 
the non-zero map $\sigma$ is an
isomorphism by Schur's lemma.

\section{Further remarks}

\subsection{Invariant sections}

\subsubsection{$G = \SL(9)$} In this case the restriction of the $E_8$-theta divisor is the unique (up
to a scalar) $\mathrm{Jac}(X)[3]$-invariant section in $H^0(\cM_X(\SL(9)), \cL_{\SL(9)})$. Here the group
$\mathrm{Jac}(X)[3]$ of $3$-torsion line bundles over X acts by tensor product on the moduli stack
of rank-$9$ vector bundles with trivial determinant. Note that we take here the linear action
of $\mathrm{Jac}(X)[3]$ on $H^0(\cM_X(\SL(9)), \cL_{\SL(9)})$ induced by the isomorphism
$\tilde{\phi}^*\cL_{E_8} = \cL_{\SL(9)}$. This leads to the natural
question (see also \cite{F2} section 6): does there exist a geometrical description of the zero-divisor of this invariant section?

\subsubsection{$G = \Spin(16)$} The $\cM_X(N)$-invariant section $\sigma \in
H^0(\cM_X(\Spin(16)), \cL_{\Spin(16)})$ can be described as follows. The center $Z$ of $\Spin(16)$
equals $\ZZ/2\ZZ \times \ZZ/2\ZZ = \{ \pm 1, \pm \gamma \}$, where $\{ \pm 1 \}$ is the
kernel $K$ of the homomorphism $\Spin(16) \rightarrow \mathrm{SO}(16)$ and $\pm \gamma$ covers the
element $-\mathrm{Id} \in \mathrm{SO}(16)$. Note that $N = \{ 1, \gamma \}$. By \cite{PR} Proposition
8.2,  the space $H^0(\cM_X(\Spin(16)), \cL_{\Spin(16)})$ admits a basis $\{ s_\kappa \}$ indexed by the
$4^g$ theta-characteristics $\kappa$ of the curve $X$ and such that the set-theoretical support
of the zero-divisor of $s_\kappa$ equals
$$ D_\kappa = \{ E \in \cM_X(\Spin(16)) \ | \ \dim H^0(X, E(\CC^{16}) \otimes \kappa ) > 0 \}.$$
Here $E(\CC^{16})$ denotes the orthogonal rank-$16$ vector bundle associated to $E$.  The group $\cM_X(N)$
identifies with the group $\mathrm{Jac}(X)[2]$ of $2$-torsion line bundles over $X$.  We can use the canonical
lift $\mathrm{Jac}(X)[2] \hookrightarrow \cG(\cL_{\Spin(16)})$ and a fixed theta-characteristic $\kappa_0$ to normalize the basis 
$\{ s_\kappa \}$ by putting $s_{\alpha \kappa}:=\alpha \cdot s_{\kappa_0}$, since any theta-characteristic $\kappa$
can be written as $\alpha\kappa_0$ for a unique $\alpha \in \mathrm{Jac}(X)[2]$.  It is then clear that 
$\sigma = \sum_\kappa s_\kappa=\sum_{\alpha \in \mathrm{Jac}(X)[2]} \alpha \cdot s_{\kappa_0}$ and that $\sigma$ does not
depend on the choice of $\kappa_0$.  As in the previous case, a geometrical description of its
zero-divisor is still missing.

\subsection{Other conformal subalgebras of $\mathfrak{e_8}$}

\subsubsection{$(\Spin(8), \Spin(8))$} We observe that $\mathfrak{p} = \mathfrak{so}(8) \oplus
\mathfrak{so}(8)$ is a {\em non-maximal} conformal subalgebra of $\mathfrak{e_8}$ with Dynkin multi-index $(1,1)$. In fact, the inclusion
$\mathfrak{p} \subset \mathfrak{e_8}$
factorizes through $\mathfrak{so}(16)$. We therefore obtain a group homomorphism $\phi: P = \Spin(8) \times
\Spin(8) \rightarrow \Spin(16) \rightarrow E_8$ with kernel $N= \ZZ/2\ZZ \times \ZZ/2\ZZ$. We note that
$N$ sits diagonally in
$Z(\Spin(8)) \times Z(\Spin(8))$. Theorem 1 and Theorem 2 also hold in this case.

\subsubsection{$(G_2,F_4)$} Using the Verlinde formula one computes that the spaces $H^0(\cM_X(G_2),\cL_{G_2})$
and $H^0(\cM_X(F_4), \cL_{F_4})$ have the same dimension, which is $\left( \frac{5 + \sqrt{5}}{2} \right)^{g-1} + \left( \frac{5 - \sqrt{5}}{2} \right)^{g-1}$. Theorem 1 says that the $E_8$-theta divisor
$\Delta$ induces
a non-zero linear map $\sigma: H^0(\cM_X(G_2),\cL_{G_2})^* \lra H^0(\cM_X(F_4), \cL_{F_4})$. Since $G_2$
and $F_4$ have no center, the argument used in the proof of Theorem 2 breaks down for this case.



\begin{thebibliography}{9999}
\bibitem[BB]{BB} A. Bais, P. Bouwknegt: A classification of subgroup truncations of the bosonic string, Nuclear
Physics B279 (1987), 561-570


\bibitem[B]{B} P. Belkale: Strange duality and the Hitchin/WZW connection, arXiv:0705.0717

\bibitem[BD]{BD} A. Borel, J. de Siebenthal: Les sous-groupes ferm\'es de rang maximum des groupes de Lie clos,
Commentarii Mathematici Helvetici, Vol. 23, No. 1 (1949), 200-221.


\bibitem[BS]{BS} A. Berenstein, R. Sjamaar: Coadjoint orbits, moment polytopes, and the Hilbert-Mumford
criterion, J. Amer. Math. Soc. 13, No. 2 (2000), 433-466

\bibitem[CG]{CG} A. Cohen, R. Griess: On finite simple subgroups of the complex Lie group of type $E_8$, The
Arcata Conference on Representations of Finite Groups (Arcata, Calif., 1986), Proceedings
of Symposia in Pure Mathematics, Volume 47, Part 2, AMS (1987)

\bibitem[D]{D} E.B. Dynkin: Semisimple subalgebras of semisimple Lie algebras, Am. Math. Soc. Transl.
(Ser. II) 6 (1957), 111-244

\bibitem[F1]{F1} G. Faltings: Theta functions on moduli spaces of $G$-bundles, J. Alg. Geom. 18 (2009), 309-369

\bibitem[F2]{F2} G. Faltings: Thetafunktionen auf Modulr\"aumen von Vektorb\"undeln, Jahresbericht der DMV, 110. Band (2008), Heft 1, 3-17

\bibitem[F3]{F3} G. Faltings: A proof for the Verlinde formula, J. Alg. Geom. 3 (1994), 347-374



\bibitem[K]{K} V. Kac: Infinite dimensional Lie algebras, Third Edition, Cambridge University Press (1990)


\bibitem[KS]{KS} V. Kac, M. Sanielevici: Decomposition of representations of exceptional affine algebras
with respect to conformal subalgebras, Physical Review D, Vol. 37, No. 8 (1988), 2231-2237

\bibitem[KW]{KW} V. Kac, M. Wakimoto: Modular and conformal invariance constraints in representation theory
of affine algebras, Advances in Math. 70 (1988), 156-234

\bibitem[KNR]{KNR} S. Kumar, M.S. Narasimhan, A. Ramanathan: Infinite Grassmannians and moduli spaces of
$G$-bundles, Math. Ann. 300 (1994), 41-75

\bibitem[LS]{LS} Y. Laszlo, C. Sorger: The line bundles on the moduli of parabolic $G$-bundles over
curves and their sections, Ann. Scient. Ec. Norm. Sup. 30 (1997), 499-525

\bibitem[MO]{MO} A. Marian, D. Oprea: A tour of theta dualities on moduli spaces of sheaves, Curves
and abelian varieties, Contemporary Mathematics, 465, AMS, Rhode Island (2008), 175-202

\bibitem[M]{M} D. Mumford: On the Equations Defining Abelian Varieties I, Invent. Math. 1 (1966),
287-354

\bibitem[NT]{NT} T. Nakanishi, A. Tsuchiya: Level-Rank Duality of WZW Models in Conformal Field Theory,
Comm. Math. Phys. 144 (1992), 351-372


\bibitem[PR]{PR} C. Pauly, S. Ramanan: A duality for Spin Verlinde spaces and Prym Theta functions, J. London
Math. Soc. 63 (2001), 513-532


\bibitem[Pa]{Pa} C. Pauly: La dualit\'e \'etrange, S\'eminaire Bourbaki No. 994, June 2008

\bibitem[Po]{Po} M. Popa: Generalized theta linear series on moduli spaces of vector bundles on curves,
arXiv:0712.3192

\bibitem[SW]{SW} A. Schellekens, N. Warner: Conformal subalgebras of Kac-Moody algebras, Physical Review D,
Vol. 34, No. 10 (1986), 3092-3096

\bibitem[So]{So} C. Sorger: On moduli of $G$-bundles on a curve for exceptional $G$, Ann. Scient. \'Ec. Norm.
Sup. 32 (1999), 127-133

\bibitem[T]{T} C. Teleman: Lie algebra cohomology and the fusion rules, Commun. Math. Phys. 173 (1995), 265-311

\bibitem[U]{U} K. Ueno: Introduction to conformal field theory with gauge symmetries, in Geometry and Physics,
Lecture Notes in Pure and Applied Mathematics 184, Marcel Dekker (1996), 603-745

\end{thebibliography}
\end{document}